\documentclass[11pt]{amsart}

\usepackage{verbatim}

\numberwithin{equation}{section}

\newtheorem{theorem}[equation]{Theorem}
\newtheorem{lemma}[equation]{Lemma}
\newtheorem{proposition}[equation]{Proposition}

\theoremstyle{definition}

\theoremstyle{remark}

\newtheorem{remark}[equation]{Remark}

\newcommand{\C}{\ensuremath{\mathbb{C}}}
\newcommand{\R}{\ensuremath{\mathbb{R}}}
\renewcommand\S{\ensuremath{{S}}}

\newcommand{\st}{\,:\,}

\newcommand\mysqrt[2]{ #2^{1/#1} }

\renewcommand\i{\mbox{\rule{2.0mm}{0.1mm}\kern0.0mm\rule{0.1mm}{3.0mm} }\,}
\newcommand\Lie[1]{{\mathcal L}_{#1}}        

\newcommand\grad{\nabla}

\renewcommand\d{\partial}

\newcommand\db{\overline{\partial}}

\newcommand\dc{d^{\C}}

\newcommand\pd[2]{\frac{\partial #1}{\partial #2}}
\newcommand\pdd[3]{\frac{\partial^2 #1}{\partial #2 \partial #3}}

\renewcommand\a{\alpha}
\newcommand\ab{\overline{\a}}

\renewcommand\b{\beta}
\newcommand\bb{\overline{\b}}

\renewcommand\c{\gamma}

\newcommand\nb{\bar{n}}

\newcommand\mapright[1]{\mbox{$\stackrel{#1}{\longrightarrow}$}}

\newcommand{\MA}{Monge-Amp\`ere}

\newcommand\Fol{\mathcal{F}}                  

\newcommand{\MC}{{\mathcal{M}}}                      

\newcommand\tFol{\widetilde{\Fol}}            

\newcommand\SM{\S M}                       

\newcommand{\tMC}{\widetilde{\MC}}

\renewcommand{\th}{\theta}
\newcommand{\tth}{\widetilde{\th}}
\newcommand{\thF}{\th_{F}}
\newcommand{\thS}{\th_{S}}

\newcommand\tc{{\tilde{\c}}}
\newcommand{\XF}{{X_{F}}}                   
\newcommand\pp{\ensuremath{\widetilde{\pi}}}       

\newcommand\tu{\tilde{u}}              

\newcommand\tp{\tilde{p}}              

\newcommand\tphi{\tilde{\phi}}         

\newcommand{\F}{F}                           

\begin{document}


\title{Singular Monge-Amp\`ere Foliations}

\author{Tom Duchamp}
\address{Department of Mathematics\\
University of Washington\\
Box 354350\\ Seattle, WA 9819-43505}
\email{duchamp@math.washington.edu}

\author{Morris Kalka}
\address{Department of Mathematics\\
Tulane University\\
New Orleans, LA 70118}
\email{kalka@math.tulane.edu}


\date{October 16, 2001}



\begin{abstract}
This paper generalizes results of Lempert and Sz\"oke on the structure of
the singular set of a solution of the homogeneous \MA\ equation on a Stein
manifold.  Their \emph{a priori} assumption that the singular set has
maximum dimension is shown to be a consequence of regularity of the
solution. In addition, their requirement that the square of the solution be
$C^3$ everywhere is replaced by a smoothness condition on the blowup of the
singular set. Under these conditions, the singular set is shown to inherit a
Finsler metric, which in the real analytic case uniquely determines the
solution of the \MA\ equation. These results are proved using techniques
from contact geometry.
\end{abstract}

\maketitle


\section{Introduction}
\label{introduction}

 The \emph{ homogeneous \MA\ equation} on the complex $n$-dimensional
complex manifold $\MC$ is the equation
\begin{equation}
\label{eqn-MA}
        (d\dc u)^n = 0 \,,
\end{equation}
where $u:\MC\to \R$ and $\dc := i(\db - \d)$. In the special case where $u$
is at least $C^{3}$ and the form $d\dc u$ has constant rank the integral
curves of $d\dc u$ foliate $\MC$ by complex submanifolds. This foliation
is called the \emph{\MA\ foliation} of $u$ and was first studied in
\cite{bedford-kalka:1977}.

An important class of solutions of \eqref{eqn-MA} is the class of
plurisubharmonic exhaustion functions for which the sets
\[
       \{ p \in \MC \colon u(p) \leq c\}
\]
are compact for all $c < \sup u$ and $d\dc u$ is a positive semidefinite
form of constant rank $n-1$. It is known (see for instance \cite[Theorem
1.1]{lempert-szoke:1991a}) that every such function must fail to be smooth
on a non-empty \emph{ singular set} $M \subset \MC$.  In this paper, we
study the extent to which the geometry of the singular set determines $u$.

\subsection*{Previous work}
Our work builds on previous results of a number of authors, particularly
those of Stoll\cite{stoll:1980}, Burns\cite{burns:1982a},
Wong\cite{wong:1982}, Patrizio-Wong\cite{patrizio-wong:1991a},
Lempert-Sz\"{o}ke \cite{lempert-szoke:1991a}, Sz\"{o}ke \cite{szoke:1991},
and Guillemin and Stenzel \cite{guillemin-stenzel:1991a}. The main result
of this paper, Theorem~\ref{thm:main}, was inspired by a question posed in
\cite{lempert-szoke:1991a}.

\medskip

In the case where $u$ has a logarithmic singularity, and $\tau=\exp(u)$ is
a smooth K\"{a}hler potential, Stoll \cite{stoll:1980} showed that $M$ is a
point and that $\MC$ is biholomorphic to either the unit ball $B^{n}\subset
\C^{n}$ or to $\C^{n}$.  Burns \cite{burns:1982a} gave a more geometric
proof, exploiting the fact that  the leaves of
the \MA\ foliation are totally geodesic with respect to the K\"{a}hler
metric.  Wong \cite{wong:1982} further explored the geometry of the \MA\
foliation, proving the following theorem.

\begin{theorem}[Wong]
\label{thm:Wong}
Let $u$ be a solution of the \MA\ equation and let $\Fol$ denote the \MA\
foliation. Let $\tau= f \circ u$ be the potential function f
for a K\"{a}hler metric on $\MC \setminus M$, where
$f$ is a smooth function satisfying the conditions $f'\circ
u>0$ and $f''\circ u>0$.   Then the leaves of $\Fol$
are totally geodesic with respect to the K\"{a}hler metric.  Moreover
if $Z$ denotes the complex vector field on $\MC \setminus M$ defined by
\[
      Z \i d\dc \tau = d\tau + i \dc \tau \,,
\]
then the leaves of $\Fol$ coincide with the orbits of the complex flow of
$Z$.  Finally, the integral curves of the real vector field $Z_{I}
=\frac{1}{2i}(Z - \overline{Z})$ are (after reparametrizing) geodesics that
intersect the level sets of $u$ at right angles.
\end{theorem}
The proof is essentially contained in \cite{wong:1982} (see also
\cite{patrizio-wong:1991a} and \cite{lempert-szoke:1991a}).

\medskip

The case where $u$ has a logarithmic singularity and $M$ is a point is an
extreme case. At the other extreme is the case where $M$ is assumed to be a
compact smooth, real $n$-dimensional submanifold of $\MC$.  Assume that $u$
is continuous on all of $\MC$ and that the singular set coincides with the
zero set of $u$.  Compactness of the set $u^{-1}([c,\infty))$ then implies
that $u$ is bounded below, we may therefore assume without loss of
generality that $u$ is non-negative and that $M$ is the zero set of $u$.
Assume, in addition, that the function $\tau = u^2$ is $C^3$ and strictly
plurisubharmonic on all of $\MC$. Then $d\dc \tau$ defines a K\"{a}hler
metric on all of $\MC$. The singular set $M$, then inherits a Riemannian
metric $g$.  The triple $(\MC,M,u)$ is  called a \emph{(Riemannian) \MA\
model}.  Patrizio and Wong \cite{patrizio-wong:1991a} studied the special
case where $M$ is a compact symmetric space. Their results were later
generalized by Lempert and Sz\"{o}ke \cite{lempert-szoke:1991a}, and
independently (when $M$ is real analytic) by Guillemin and Stenzel
\cite{guillemin-stenzel:1991a}, to the case where $M$ is an arbitrary
compact Riemannian manifold. The results in \cite{lempert-szoke:1991a} were
extended by Sz\"{o}ke in  \cite{szoke:1991} and
\cite{szoke:1995}. The main results of
\cite{lempert-szoke:1991a} and \cite{szoke:1991} are summarized in
the following theorem:

\begin{theorem}[Lempert-Sz\"{o}ke \cite{lempert-szoke:1991a} and  Sz\"{o}ke \cite{szoke:1991}]
\label{thm:1}
(a)
 Let $(\MC,M,u)$ be a \MA\ model. Then the set of curves given by the
intersection of $M$ with the leaves of the \MA\ foliation is precisely the
set of geodesics of $M$ with respect to the induced metric $g$.

(b) Every compact real analytic Riemannian manifold arises from a \MA\
model.

(c) Let $(\MC,M,u)$ and $(\MC',M',u')$ be two \MA\ models. Suppose that
$(M,g)$ and $(M',g')$ are isometric and $\sup u = \sup u'$ then there is a
biholomorphic map $\Phi:\MC \to \MC'$ such that $u = u' \circ \Phi$.

\end{theorem}

In a related paper \cite{lempert:1993}, Lempert showed that the Riemannian
manifold $M$, metric $g$, and exhaustion function $u^2$, associated to a
\MA\ model are all real analytic. And in \cite{szoke:1995} Sz\"{o}ke proved
a further generalization of part (c) of Theorem~\ref{thm:1}.

\subsection*{Results}
Our goal is to understand the structure of the singular set of a solution
of the \MA\ equation under weakened smoothness assumptions on $u$ as well
as weakened assumptions on the topology of $M$.  Throughout this paper
$\MC$ denotes a complex $n$-dimensional Stein manifold. We remark that the
assumption that $\MC$ is Stein is made to avoid cases such as $\MC=X\times
Y$ with $X$ Stein and $Y$ compact.  When we say that $u$ \emph{is a
solution of the homogeneous \MA\ equation\/} we always mean that $u$ is an
everywhere continuous, non-negative, plurisubharmonic exhaustion of $\MC$
that is a solution of the equation
\[
       (d\dc u)^n = 0 \,, \quad (d\dc u)^{n-1} \neq 0
\]
on the set $\MC \setminus M$, where $M$, the zero set of $u$, is assumed to
be a smooth compact submanifold. We also assume the $u$ is smooth on the
complement of $M$\footnote{Most of our results apply to the case where $u$
is only of class $C^3$ , the minimum smoothness assumption needed to make
our geometrical constructions.}

Additional smoothness assumptions on $u$ and $M$ are made in both
\cite{lempert-szoke:1991a} and \cite{guillemin-stenzel:1991a}, where $u^2$
is assumed to be a smooth K\"{a}hler potential on all of $\MC$ and $M$ is
assumed to have real dimension $n$.
In Section~\ref{sec:singular}, we show that both of these assumptions can
be weakened, Theorem~\ref{thm:smoothness} shows that the assumption that
$u^2$ is a smooth K\"{a}hler potential implies a regularity condition for
$u$ on the normal-blowup of $M$ (see below). And Theorem~\ref{thm:singular}
shows that the regularity condition implies that the singular set is an
$n$-dimensional, totally real submanifold.

Our  regularity condition is expressed in term of
the \emph{normal blowup of $M$ in $\MC$}, which is a smooth manifold with boundary
$\tMC$, together with a smooth map
\[
    \pp:\tMC \to \MC
\]
such that (i) the space $\SM = \pp^{-1}(M)$ is diffeomorphic to the normal
sphere bundle of $M$ in $\MC$, and (ii) $\pp: \tMC \setminus \SM \to \MC
\setminus M$ is a diffeomorphism.  Section~\ref{sec:blowup} contains a more
detailed description.  Let $\tu = u \circ \pp$ and $\tth = \pp^{*}\dc u$.
We assume that $\tu$ is a smooth function.%
\footnote{Most of our
computations require that $\tu$ only be $C^3$.}
And we replace the
assumption that $u^2$ is strictly plurisubharmonic by the assumption that
the 1-form $\th$ extends smoothly to all of $\tMC$ and satisfies the
non-degeneracy condition
\[
         d\tu \wedge \tth \wedge (d\tth)^{n-1} \neq 0\,.
\]
When $u$ satisfies these conditions we say that $u$ is \emph{regular on the
normal blowup of $M$}.

\medskip

Theorem~\ref{thm:1} generalizes to the case where $u$ is regular on the
normal blowup.  We prove that regularity implies that  the pull-back to $\SM$ of the
form $\tth$ is a contact form on $\SM$.  Let $Q$ denote the normal bundle
of $M$. Then, by Theorem~\ref{thm:singular}, the bundle map
\[
     TM \mapright{J} \left.T\MC \right|_{M} \mapright{\pi_Q} Q\,,
\]
where $J$ is the complex structure tensor of $\MC$ and $\pi_Q$ is
projection onto the normal bundle, is an isomorphism. Therefore, $\SM$ can
be identified with the projective tangent bundle of $M$. We show in
Section~\ref{sec:finsler} that this form defines a Finsler metric $\F$ on
$M$. Following the terminology of \cite{patrizio-wong:1991a}, we say that
the triple $(\MC,M,u)$ a \emph{regular \MA\ model} for the Finsler metric
$\F$. Our main result is the following theorem.

\begin{theorem}
\label{thm:main}
Let $\MC$ be a Stein manifold and let $u\geq 0$ be a solution of the \MA\
equation $(d\dc u)^n = 0$, $(d\dc u)^{n-1} \neq 0$ on $\MC\setminus M$,
where $M = \{ u=0\}$. Finally assume that $M$ is a compact, smooth
submanifold.

(a) If $u$ is regular on the normal blowup of $M$, then $(\MC,M,u)$ is a
regular \MA\ model for a Finsler metric $(M,\F)$. The leaves of the \MA\
foliation intersect $M$ along geodesics.

(b) Every real analytic Finsler metric on $M$ arises from a regular \MA\
model.

(c) Let $(\MC,M,u)$ and $(\MC' ,M',u')$ be two real analytic \MA\ models
for the real analytic Finsler metrics $(M,\F)$ and $(M',\F')$,
respectively.  Then there is a biholomorphic map $\Phi:\MC \to \MC'$,
defined in a neighborhood of $M$ such that $u = u' \circ \Phi$ if and only
if $(M,\F)$ and $(M,\F')$ are isometric.
\end{theorem}

The paper is organized as follows. In Section~\ref{sec:contact} we show
that the space $\MC \setminus M$ is diffeomorphic to the product of a
contact manifold and an open interval. In Section~\ref{sec:blowup} we
define the normal blowup. In Section~\ref{sec:singular}, we give a precise
definition of regularity on the normal blowup, extend the contact structure
in Section~\ref{sec:contact} to the the blowup, and prove that regularity
implies total reality of $M$. In Section~\ref{sec:finsler}, we review
Finsler geometry and give the proof that $M$ inherits a Finsler metric. In
Section~\ref{sec:models}, we complete the proof of Theorem~\ref{thm:main}.

\medskip

\begin{remark}
 In \cite{lempert-szoke:1991a}, Lempert and Sz\"{o}ke conjectured that the
case when $u^2$ is not smooth could be studied by replacing the Riemannian
metric $g$ by a Finsler metric. Theorem~\ref{thm:main} confirms this
conjecture.
\end{remark}

The authors would like to thank the referee for pointing out a
problem in the proof of part c) of  \ref{thm:main}


\section{Contact geometry away from the singular set}
\label{sec:contact}

Our  goal is to understand the relation between
solutions of the \MA\ equation on a  Stein manifold
and the geometry of its singular set.

We assume that $\MC$ is a complex $n$-dimensional Stein manifold and that
$u$ is an everywhere continuous, non-negative, solution of the homogeneous
\MA\ equation
\[
         (d\dc u)^n = 0
\]
whose zero set $M = \{ u = 0\}$ is a smooth, compact submanifold.
By restricting to a neighborhood of $M$ if necessary, we assume that
 $u$ is bounded above by $R>0$ and that
$u^{-1}(r)$ is compact for all $0 \leq r < R$.
Finally, we assume that $\tau = u^2/2$ is strictly plurisubharmonic on $\MC
\setminus M$.\footnote{We gain no further generality by
replacing $\tau=u^2/2$, by a more general function $\tau = f(u)$, with
$f'(u), f''(u)>0$ for $u>0$: For  $(d\dc\tau)^n = n (f')^{n-1} f'')
 du \wedge \dc u \wedge (d\dc u)^{n-1}$, which  is a positive multiple
of $(d\dc (u^2/2))^n$.}

Because $\tau$ is strictly plurisubharmonic, the two form $d\dc\tau$ has
rank $n$ away from $M$, and the computation
\[
(d\dc \tau)^n =
\left(  du \wedge \dc u +  u\, d\dc u \right)^n
=
 \left( du \wedge \theta + u\, d\theta \right)^n
= n \,u^{n-1}\, du \wedge \th \wedge (d \th)^{n-1} \,.
\]
shows that $u$ satisfies the non-degeneracy condition
\begin{equation}
\label{eq:contact-1}
          d u \wedge \dc u\wedge (d\dc u)^{n-1} \neq 0\,.
\end{equation}
It follows that the level set $M_{\epsilon} = \{ u = \epsilon \}$ is a
smooth, contact manifold for all $\epsilon$ between $0$ and $R$. The
contact form on $M_{\epsilon}$ is the pull-back of the one-form
\[
             \th = \dc u \,.
\]
Because $u$ has no critical points, the level sets $M_{\epsilon}$ are all
diffeomorphic. Indeed, they are isomorphic as calibrated, contact
manifolds. (A \emph{calibrated contact manifold} is a contact manifold with
a distinguished contact form.) To see this, let $Y$ be the vector field
characterized by the conditions
\begin{equation}
\label{eq:contact-2}
           Y \i \theta = 0 \,, Y \i du = 1\,,\text{ and } Y \i d\theta = 0 \,.
\end{equation}
Let $\mu_{t}$ denote the flow of $Y$. Because $du(Y)=1$, $\mu_t$ maps level
sets of $u$ to level sets, and, therefore, defines a diffeomorphism
\begin{equation}
\label{eq:contact-3}
      \mu:   M_{R/2} \times (0, R) \to \MC \setminus M \,:\, (p,t)
         \mapsto \mu_{t-R/2}(p) \,
\end{equation}
satisfying the identity
\begin{equation}
\label{eq:contact-4}
     u \circ(\mu(p,t) ) = t \,.
\end{equation}
The  computation of the Lie derivative
\[
       \Lie{Y} \th = d (Y \i \th) + Y \i d\th = 0
\]
then shows that $\mu_t$ restricts to a contact diffeomorphism
\[
      \mu_{\epsilon_2-\epsilon_1}: M_{\epsilon_1} \to M_{\epsilon_2}\,.
\]
between each pair of level sets.
\medskip

\begin{remark}
\label{rem:contact}
 The form $\th$ satisfies an even stronger condition: The identities
$Y\i\th =0$ and $ \Lie{Y}\th = 0$ together imply that $\th$ descends to a
contact form on the orbit space $(\MC \setminus M)/Y \simeq M_{R/2}$. This
implies that $\MC \setminus M$ has the structure of the product of a
contact manifold with the interval $(0, R)$ and that the pull-back
$\mu^{*}\th$ extends smoothly to $M_{R/2} \times [0,R)$ with
$\mu^{*}(\th\wedge(d\th)^{n-1} \neq 0$ everywhere.
\end{remark}

\section{The normal blowup of $M$}
\label{sec:blowup}

Because $u$ is continuous, the level sets $M_{\epsilon}$ approach $M$ as
$\epsilon$ approaches $0$. When $u$ is sufficiently well-behaved, the
contact structures on the level sets converge to a limiting contact
structure on the projective normal bundle of $M$.
The normal blowup of $M$, defined below, is our main tool for formalizing this behavior.

We first discuss the simpler case of the blowup  of the origin in $\R^m$.  In
this context, blowing up is just the transformation to spherical
coordinates. Spherical coordinates, which we formalize by the
\emph{blowdown map}
\[
  \pp: \widetilde{\R}^q := S^{q-1} \times [0,\infty) \to \R^m: (v,r)
           \mapsto   r \cdot v \,,
\]
where $S^{q-1}$ is the unit sphere in $\R^q$.  The preimage
 $\pp^{-1}(0)$ is called the \emph{blowup of the origin}. Notice that any smooth
curve satisfying the conditions $\c(t)$, $t \geq 0$ with $\c(0)= 0$, $\c'(0)
\neq 0$, and $\c(t) \neq 0$ for $t>0$, has a unique lift to a smooth curve
 on the blowup defined by
\[
      \tc(t) =  \begin{cases}
                     \left(  \frac{\c(t)}{\|\c(t)\|}, \|\c(t)\| \right) &
                     \text{for $t>0$,}\\
                    \left(  \frac{\c'(0)}{\|\c'(0)\|}, 0 \right) &
                        \text{ for  $t=0$.}
                \end{cases}
\]
\begin{remark}
\label{rem:blowup}
We want to emphasize the following three obvious properties of the lift:
\begin{enumerate}
\item[ (i)] $\tc(t)$ intersects the boundary of $\widetilde{\R}^m$
transversely;
\item[(ii)]$\tc(0)$ depends only on the oriented ray generated by $\c'(0)$;
\item[(iii)]$\tc(0) = \c'(0)/ \|\c'(0)\|$.
\end{enumerate}
\end{remark}

Roughly speaking, the normal blowup of a submanifold $M$ is obtained by
replacing each point of $M$ by the blowup of the origin of the vector space
of normal vectors to $M$ in $\MC$.  We now present a more formal
description.

Consider first the case where $V$ is a $q$-dimensional vector space and $M$
is the origin. Let $V_0$ be the set of non-zero vectors, and let $\S V$
denote the space of oriented rays through the origin. We call $\S V$ the
\emph{(oriented) projectivization} of $V$.  The \emph{blowup} of $V$ at the
origin is the subspace
\[
      \widetilde{V} =   \{ ([v], r \cdot v) \in \S V \times V \st v\in V, v
      \neq 0,\, r \in [0, \infty) \} \,,
\]
where $[v]$ denotes the oriented ray defined by the non-zero vector $v \in
V$.  This definition generalizes fiber wise to a vector bundle $E$ in the standard
way. In this case $\S E$ denotes the oriented projective bundle of $E$ and
$\widetilde{E}$ denotes the blowup of the set of zero vectors of $E$. There is
a natural \emph{blowdown} map $p: \widetilde{E} \to E$.
It is easy to check that if $E$ is equipped with a norm, we can identify
$\S E$ with  the set of unit length vectors, and the map
\[
    SE \times [0,\infty) \to \widetilde{E} \st (v,r)
          \mapsto ([v], r \cdot v)
\]
is a diffeomorphism. In particular, $\S E$ is a sphere bundle over $M$.
The  canonical map
\[
    \pp:  \widetilde{E} \to E
\]
sending $([v],v)$ to $v$ and $([v],0)$ to the zero vector is called the
\emph{blowdown} map.

\medskip

The normal blowup of a submanifold is the non-linear version of the blowup
of the zero section of a vector bundle.  We give  two equivalent
constructions here.  The first highlights the role of the
normal bundle and uses the exponential map of an auxiliary metric, the
second is based on local coordinate charts and does not rely on an explicit
choice of metric. The proof that these constructions are equivalent is an
exercise in differential geometry, which we leave to the reader.

Let $Q$ denote the normal bundle of $M$. Then there is a short exact
sequence of vector bundles
\[
        0 \to TM \longrightarrow \left. T\MC \right|_{M} \mapright{\pi_Q} Q
        \to 0 \,.
\]
A choice of a Riemannian metric on $\MC$ gives a splitting, under which $Q$
can be identified with the orthogonal complement of $TM$ in $T\MC$.  The
exponential map defines a diffeomorphism between an $\epsilon$-neighborhood
of the zero-section of $Q$ and a neighborhood of $M$ in $\MC$.  Let
$\widetilde{B}_{\epsilon} \subset \widetilde{Q}$ be a neighborhood of the
blowup of the zero-section of $Q$.  The \emph{normal blowup} of $\MC$ along
$M$ is the manifold $\tMC$ obtained by identifying points in the manifold
$\MC \setminus M$ with points in $\widetilde{B}_{\epsilon}$ by the
exponential map. Let
\[
       \pp: \tMC \to \MC
\]
be the blowdown map, defined in the obvious way.  Notice that
$\tMC$ is a smooth manifold whose  boundary is the subspace $\SM =\pp^{-1}(M)$.
By definition, $\SM = \S Q$.  Observe also that the distance to $\SM$ is comparable to
the distance to $M$ with respect to the Riemannian metric on $\MC$.
We call the submanifold  $\SM$ the \emph{normal blowup of $M$} (or less
formally, the \emph{blowup of $M$}).

Our second  construction of the blowup begins with a collection $U_{\a}$ of open
subsets of $\MC$ whose union contains $M$, together with a collection of
coordinates charts
\[
    \phi_{\a} : U_{\a} \to V_{\a} \times B_{\epsilon}^q   \st p \mapsto (x,y),
\]
which satisfy the compatibility condition $M \cap U_{\a} = V_{\a}
\times \{0\}$, where $V_{\a}$ is an open subset of $\R^n$ and
$B^q_{\epsilon}$ denotes the  ball of radius $\epsilon$ in $\R^q$ centered
at the origin. The transition functions are maps of the form
\[
\phi_{\a,\b}= \phi_{\b} \circ \phi^{-1}_{\a}
\st V_{\a,\b} \times B_{\epsilon}^{q} \to
\st V_{\b,\a} \times B_{\epsilon}^{q}
\st  (x,y) \mapsto  (X(x,y),Y(x,y))
\]
where $V_{\a,\b} = \phi_{\a}\left(M \cap U_{\a} \cap
U_{\b}\right)$. By virtue of the compatibility condition, the
$y$-component of the transition functions can be written in the form
\begin{equation}
\label{eq:blowup-1}
      Y^{k}(x,y) =  a^{k}_{i}(x) y^{i} +  R^{k}_{i,j}(x,y) y^i y^j
\end{equation}
where $A = ( a^i_j (x) )$ is a smooth family of invertible $q\times q$
matrices and $R^{k}_{i,j}(x,y)$ are smooth functions, the indices $i,j,k$
ranging between $1$ and $q$ with the summation convention in force.  Thus,
for $t \geq 0 $ sufficiently small, the transition functions induce maps
\begin{equation}
\label{eq:blowup-2}
   \widetilde{\phi}_{\a,\b}:
 V_{\a,\b} \times  S^{q-1} \times [0,\epsilon)
\to
 V_{\b,\a} \times  S^{q-1} \times [0,\epsilon)\,,
\end{equation}
defined by the formula
\begin{equation}
\label{eq:blowup-3}
\widetilde{\phi}_{\a,\b}(x,v, r) =
\begin{cases}
 \left( X(x, r v), \frac{ Y(x,r v) }{ \| Y(x, r v) \|  }, \| Y(x,r v) \| \right)
& r>0\,,\\
 \left(  X(x, 0), \frac{ A \cdot v }{ \|  A \cdot v \| }, 0 \right)
& r = 0   \,.
\end{cases}
\end{equation}
A straightforward computation shows that these functions satisfy the
cocycle condition
\[
       \widetilde{\phi}_{\b,\c} \circ\widetilde{\phi}_{\a,\b}
        = \widetilde{\phi}_{\a,\c} \,.
\]
Let $\sim$ denote the equivalence relation  on the disjoint union
$
(\MC \setminus M)\, \overset{\cdot}{\cup}\,
\overset{\cdot}{\bigcup_{\a}}\, \left( V_{\a} \times S^{q-1}  \times [0,\epsilon)\right)
$
generated by the relations $(x,v,r)\sim\widetilde{\phi}_{\a,\b}(x,v,r)$ and
$p \sim\phi_{\a}\circ \pp(p)$. The \emph{normal blowup of $M$ in
$\MC$} is defined to be the quotient space
\[
(\MC \setminus M)\, \overset{\cdot}{\cup}\,
\overset{\cdot}{\bigcup_{\a}}\, \left( V_{\a} \times S^{q-1}  \times [0,\epsilon)\right)
\,/\sim \,.
\]
The cocycle condition guarantees that $\widetilde{\MC}$ is a smooth
\mbox{$(n+q)$-dimensional} manifold with boundary diffeomorphic to $\S Q$;
and the blowdown map $\pp$ is smooth by construction.

\medskip

The verification that  two definitions of normal blowup are equivalent is
an elementary exercise in differential geometry, which we leave to the
reader.

\subsection{Blowup coordinates.}
We will often have to work in local homogeneous coordinates centered at an
oriented normal ray in $\SM$. More specifically, we shall choose a local
coordinate chart $\phi:U \to V \times B_{\epsilon}^q$ with local coordinate
functions
\[
    (x,y) = \left(x^1,\dots,x^n, y^1,\dots,y^{q} \right)
\]
on $\MC$ such that $M$ intersects $U$ in the set $\{ y = 0\}$. The collection
of points of $\tMC$ over $U$ is then a set of the form
\[
         \left( x, \frac{y}{\|y\|}, \|y\| \right) \in V \times S^{q-1}
         \times [0,\epsilon) \,.
\]
We shall choose $\phi$ so that the ray of interest is defined by $y =
(0,\dots,0,1)$. The map
\[
          \left( x, \frac{y}{\| y\|},\|y\| \right) \mapsto
          (x, p, r) = (x, \frac{y^1}{y^q},\dots, \frac{y^{q-1}}{y^q}, y^q)
             \in V \times R^{q-1} \times [0,\infty)
\]
is clearly a coordinate chart for $\tMC$ centered at the ray. We shall refer to such
coordinates as \emph{blowup coordinates}.  In blowup coordinates, the
blowdown map assumes the form
\begin{equation}
\label{eq:blowup-4}
  \pp: (x,p,r) = (x,p^1,\dots,p^{q-1},r)
\mapsto (x, y ) = (x, (r\, p^1, \dots, r \, p^{q-1},r)) \,.
\end{equation}

\bigskip

The following lemma   summarizes some of the elementary  properties of the
blowup that we need. It is an obvious extension of
Remark~\ref{rem:blowup}. The proof is an elementary exercise, which we leave to the reader.

\begin{lemma}
\label{lem:blowup}
Let $\c(t)$, $ t\geq 0$,  be a smooth curve in $\MC$ intersecting $M$ transversely
at $t=0$, with $\c(t)\notin M$ for $t>0$.

\begin{enumerate}
\item[(i)] Then  $\c(t)$ has a unique lift to
a smooth curve $\widetilde{\c}(t)$ in $\tMC$ defined by
 letting $\tc(0) \in \S Q$ be the  oriented ray generated by
$\pi_Q(\c'(0))$.

\item[(ii)] Let $f$ be a smooth function on
$\tMC$ that vanishes on $\SM$, then the quantity
$ df(\widetilde{\c}'(0) )$ depends only on  $\pi_Q \c'(0)$.

\item[(iii)] Let $Y$ be a vector based at a point
$p\in \SM = \S Q$ such that $ \pi_{Q} \left( \pp{}_* Y \right) \neq 0$.
Then the ray defined by the normal vector
$\pi_{Q} \left( \pp{}_* Y \right)$ is $p$,
itself.
\end{enumerate}
\end{lemma}


\section{The structure of the singular set}
\label{sec:singular}

In this section, we give a regularity condition on $u$ that generalizes the
one given in \cite{lempert-szoke:1991a} and explore some of its
implications.  Set $\tu = \pp^{*}u$ and $\tth = \pp^{*}\th$.  We say that
$u$ is \emph{regular on the normal blowup of $M$} (or more simply
\emph{regular on the blowup}) if and only if is satisfies the following two
conditions:
\begin{enumerate}
\item[(i)] $\tu$ and $\tth$ extend smoothly to all of $\tMC$,
\item[(ii)] the form $d\tu \wedge \tth \wedge (d\tth)^n$ is non-vanishing on
all of $\tMC$.
\end{enumerate}

\medskip
The next proposition roughly states that  regularity on the blowup is equivalent to
the condition that $\tMC$ be the product of a contact manifold with an
interval. This is the main geometric fact underlying all of our results.

\begin{proposition}
\label{prop:smoothiness-1}
The diffeomorphism  of Equation~\eqref{eq:contact-3} extends to a diffeomorphism
\[
       \widetilde{\mu}: M_{\epsilon} \times [0,R) \to \tMC
\]
if and only if $u$ is regular on the blowup.
\end{proposition}

\begin{proof}
Assume that $\mu$ extends to a diffeomorphism $\widetilde{\mu}$ as above;
then, by virtue of Equation~\eqref{eq:contact-4}, $\tu = \pi_{2}
\circ\widetilde{\mu}$, where $\pi_2(p,t) = t$.  Because $\pi_2$ is smooth
and $d\pi_2 = dt$, the function $\tu$ is smooth and $d\tu$ never
vanishes. Recall from Remark~\ref{rem:contact}, that the form $\mu^*(\th)$
extends smoothly to all of $M_{\epsilon} \times [0,R)$ and restricts to a
contact form on $M_{\epsilon} \times \{0\}$. This implies that $\tth$ is smooth
on all of $\tMC$ and restricts to a contact form on $\SM$.

Conversely, suppose that $\tu$ is regular on the blowup and that the form
$\tth$ is extends smoothly to all of $\tMC$ and restricts to a contact form
on $\SM$.  Then because $d\tu$ is non-vanishing, $M_{R/2}$ is
diffeomorphic to $\SM$. It also follows that the construction of the vector field $Y$
given in Section~\ref{sec:contact} extends to define a vector field
$\widetilde{Y}_u$  on all of $\tMC$. Since $\widetilde{Y}_u$ is
transverse to $\SM$, the map
\[
  \widetilde{\mu}: \SM \times [0,R) \to \tMC
\st
(p,t) \mapsto \widetilde{\mu}_t(p)
\]
where $\widetilde{\mu}_t$ is the flow of $\widetilde{Y}_u$ is a
diffeomorphism. By uniqueness of integral curves, $\widetilde{\mu}$
agrees with $\mu$ on the interior of $\tMC$.
\end{proof}

\begin{remark}
A result very much like this appears in the paper of Burns \cite{burns:1982a}.
\end{remark}

\medskip

Recall that the Theorems of Stoll \cite{stoll:1980} and Lempert-Sz\"{o}ke
\cite{lempert-szoke:1991a} concern the structure of the singular set of $u$
in the extreme cases where its dimension is either $0$ or $n$. Our next
result shows that under mild regularity conditions on $u$, no other
dimensions are possible.

\begin{theorem}
\label{thm:singular}
Suppose that $u$ is a solution of the \MA\ equation that is regular on the
 normal blowup of $M$.  Then $M$ is an $n$-dimensional, totally real
 submanifold of $\MC$.
\end{theorem}

Our proof proceeds by studying the lift of the \MA\ foliation $\Fol$ to
$\tMC$.  Assume that $u$ is regular on the blowup. Then by
Proposition~\ref{prop:smoothiness-1}, the closed form $d\tth$ has rank
$n-1$ everywhere on $\tMC$, as does its restriction to $\SM$, the boundary
of $\tMC$. Consequently, $\Fol$ lifts to a non-singular foliation $\tFol$
of $\tMC$ by (real) surfaces, and the leaves of $\tFol$ intersect $\SM$
transversely in curves.

Each leaf of $\Fol$ has a holomorphic parameterization expressed in
terms of the complex flow of the complex vector field
\begin{equation}
\label{eq:singular-1}
        Z =    X   + i \, Y\,,
\end{equation}
where $X$ and $Y$ are real vector fields on $\tMC$ characterized by the conditions
\begin{equation}
\label{eq:singular-2}
           X \i \tth = -1 \,, X \i d\tu = 0\,, X \i d\tth = 0 \,
\text{ and }
           Y \i \tth = 0 \,, Y \i d\tu = 1 \,, Y \i d\tth = 0 \,.
\end{equation}
Notice that  $Y$ is the extension to all of $\tMC$ of
the vector field defined in Equation~\eqref{eq:contact-2}.
Let $\nu_t$ and $\mu_t$ be the flows of $X$ and $Y$, respectively. For each
point $\tp \in \SM$, consider the map
\begin{equation}
\label{eq:singular-3a}
    \tphi_{\tp}: H \to \tMC
       \st \zeta = s + i r \mapsto \nu_{s}\circ \mu_{r}(\tp) \,,
\end{equation}
where $H = \{ s + i r \in \C \st 0 \leq r < R \}$, and set
\begin{equation}
\label{eq:singular-3b}
        \phi_{\tp} = \pp \circ \tphi_{\tp}: H \to \MC \,.
\end{equation}

\begin{lemma}
\label{lem:singular}
For each $\tp \in \SM$, the map $\tphi_{\tp}$ is well-defined and the map
$\phi_{\tp}$ is holomorphic and
non-singular at all points of $H$. The collection of images of $\phi_{\tp}$
as $\tp$ ranges over all of $\SM$ spans the \MA\ foliation.
Finally, the leaf of $\Fol$ defined by $\phi_{\tp}$ intersects $M$ along
the non-singular curve
\[
           s \mapsto   \phi_{\tp}( s) \,, s \in \R \,.
\]
\end{lemma}

\begin{proof}
We claim that $X = JY$ away from $\SM$ and that $X$ and
$Y$ commute everywhere.
To verify the first condition, recall that $\tth$ is the extension of $\th
= \dc u$ to $\tMC$, and that $\dc = d \circ J$, where $J$ is the complex
structure tensor; the identity $Y = JX$ follows from the definition
\eqref{eq:singular-2}.
To see that the vector fields $X$ and $Y$ commute, recall that because
$d\tu \wedge \tth \wedge (d\tth)^{n-1}$ is a volume form on $\tMC$, we need
only prove that the Lie bracket $[X,Y]$ is in the kernel of each of the
forms $d\tu$, $\tth$, and $d\tth$. But
\begin{align*}
 0 &= d^2\tu(X,Y) = X d\tu(Y) - Y d\tu(X) - d\tu([X,Y]) = -d\tu([X,Y])\\
\intertext{and}
0 &= d\tth(X,Y) = X\tth(Y) - Y\tth(X) - \tth([X,Y]) = - \tth([X,Y]) \,,
\end{align*}
and, because the kernel of $d\tth$ is an involutive distribution and $X$
and $Y$ are both in the kernel, so is $[X,Y]$.
Because $\MC = u^{-1}( [0, R) )$ and the level sets of $u$ are all compact,
$\tphi$ is well defined on all of $H$.  Because $Y = JX$ on $\MC \setminus
M$, $\phi$ is a holomorphic curve in $\MC$.

By construction, the vectors $\pp{}_{*}(X)$ and $\pp{}_{*}Y(p)$ are non
vanishing for all $p \in \tMC$. Consequently, $\phi$ is a non-singular
parameterization of a leaf of the \MA\ foliation.

To see  that every leaf of $\Fol$ is contained in the image of
$\phi_{\tp}$ for some $\tp \in \SM$, choose a point in $ p \in
\MC \setminus M = \tMC \setminus \SM$. Then $p = \mu_{r}(\tp)$ for a unique
point $\tp \in \SM$. Hence, the leaf of $\Fol$ through $p$ is contained in
the image of $\phi_{\tp}$.

Finally, to verify that the leaves of $\Fol$ intersect $M$ along  non-singular
curves of the form $s\mapsto \phi_{\tp}(s)$, recall that $X(\tu)=0$. This
shows that the curve is contained in $M$. Moreover, by construction,
\[
     \phi_{\tp}'(s) = J \left(\pp{}_{*}Y_{\tphi(s)}\right) \neq 0 \,,
\]
showing that the curve is non-singular.
\end{proof}

\begin{proof}[Proof of Theorem~\ref{thm:singular}]
Assume that $u$ is regular on the normal blowup.
First observe that the flow of $Y$ induces a continuous deformation retract
$\rho:\MC \to M$
defined as follows
\[
  \rho(p) = \begin{cases}
                p & \text{for $p\in M$ }\\
            \pp \circ  \mu_{-u(p)}\left(\pp^{-1}(p)\right) &\text{for $p \in \MC \setminus
                M$.}
            \end{cases}
\]
Hence, $M$ and $\MC$ have the same homotopy type. By the theorem of
Andreotti-Frankel~\cite{andreotti-frankel:1959}, the Stein manifold $\MC$
has the homotopy type of an $n$-dimensional cell complex. Consequently,
$M$ can have dimension at most $n$.

Let $TM$ denote the tangent bundle of $M$, and let $J:T\MC \to T\MC$ denote
the complex structure tensor of $\MC$. We claim that the composition
\begin{equation}
\label{eq:singular-4}
          TM \mapright{J} \left.T\MC \right|_{M} \mapright{\pi_Q} Q
\end{equation}
is a surjective map onto the normal bundle of $M$ in $\MC$. Because the
dimension of $M$ is at most $n$, this claim
implies, that the map \eqref{eq:singular-4} is an isomorphism of vector
spaces, hence, that $M$ is totally real.

To prove that the map \eqref{eq:singular-4} is surjective, first choose a point
$p\in M$ and a non-zero vector $v \in Q_p$. We need only show that a
multiple of $v$ is in the image of this map. But the vector $v$ defines an
oriented ray, which by definition of $\SM$ is a point $\tp \in \SM$ with
$\pp(\tp) = p$.  By Lemma~\ref{lem:blowup}(iii), the oriented rays defined by
$\pi_Q \left( J \pp{}_{*}X_{\tp} \right)$ and  $v$ coincide.
\end{proof}

\begin{remark}
 A theorem of Harvey and Wells \cite{harvey-wells:1973} states that the zero
set of a non-negative, strictly plurisubharmonic function is locally
contained in a totally real submanifold. Because $u^2$ may not be smooth on
$\MC$, the theorem  does not apply.
\end{remark}

\medskip

Our next theorem shows that Proposition~\ref{prop:smoothiness-1} is,
indeed, a generalization of the requirement in \cite{lempert-szoke:1991a}
that $u^2$ be a smooth K\"{a}hler potential.  More generally, one could
assume that the function $\tau = f(u)$ is a smooth potential function for a
K\"{a}hler metric.  The next theorem shows that all such conditions imply
that $u$ is regular on the blowup.

\begin{theorem}
\label{thm:smoothness}
Let $u \geq 0$ be a solution of the \MA\ equation on $\MC$. Assume that the singular
set $M = \{ u=0\}$ a smooth submanifold.  Suppose that $\tau= f\circ u \geq
0$ is a smooth, strictly plurisubharmonic exhaustion function for $\MC$,
where $f$ is a real analytic function with  $f(0) = 0$ and with $f'(u)$ and
$f''(u)$ both positive for $u>0$. Then $u$ is regular
on the blowup and $M$ is a totally real submanifold of maximum dimension.
\end{theorem}

\begin{proof}
We first claim that the form $\tth$ extends smoothly to all of $\tMC$.  To
see this, give $\MC$ the K\"{a}hler metric defined by the K\"{a}hler
potential $\tau$.  One easily verifies that the vector field $Z = X + i Y$
defined in \eqref{eq:singular-1} satisfies the identity
\[
       Z \i d\dc\tau = \frac{f''(u)}{f'(u)} \, (d\tau + i \dc\tau)
\]
on $\MC \setminus M$.  Therefore by Theorem~\ref{thm:Wong}, the gradient
vector field $\grad \tau$ is a scalar multiple of $Y$, and each integral
curve of $Y$ is contained in a geodesic of $\MC$ that intersects the level
sets of $\tau$ orthogonally.

Consequently, these geodesics lift to the blowup and intersect the boundary
of $\tMC$ transversely. The union of all of these curves forms a one
dimensional foliation of $\tMC$ with tranversal intersection with the
boundary of $\tMC$. Moreover, the leaves of this foliation are (by
construction) the closures of the integral curves of $Y$.  The identities
$Y\i\tth = \Lie{Y}\tth =0$ then show that the form $\tth = \pp^*\dc u$
extends smoothly to all of $\tMC$ and is non-vanishing at all points of
$\SM$.

Let $\thS$ denote the pullback of $\tth$ to the boundary $\SM \subset
\tMC$. The non-degeneracy condition $\tth \wedge (d\tth)^{n-1} \neq 0$,
implies that $\thS$ is a contact form on $\SM$.  Therefore, to conclude the
proof of regularity, we need only show that $\tu$ is smooth on all of
$\tMC$ and that $d\tu$ is non-vanishing near $\SM$. We do this obtaining
explicit formulas for $\dc u$ and $\tu$ in blowup coordinates adapted to
the complex structure on $\MC$.

By a theorem of Harvey and
Wells \cite{harvey-wells:1973} (see also \cite{lempert-szoke:1991a}), $M$,
the zero set of a smooth strictly plurisubharmonic function, is totally
real. Let $m \leq n$ be the dimension of $M$, and let $q = n-m$, and let
the indices $j$ and $a$ range between $1$ and $m$ and $1$ and $q$,
respectively.

We choose holomorphic coordinates
\[
          \MC \supset  U \to \C^{m+q} \st p \mapsto (z^1,\dots,z^{n} )
\]
with $z = x + i y$ and a smooth function $H:\R^m \to \R^m \times \C^q$ such that
\[
    M \cap U = \{ z \in \C^{m+q}
       \st (y^1,\dots,y^m,z^{m+1},\dots,z^{m+q})
        =   H(x^1,\dots,x^m)  \}  \,.
\]
Because $M$ is totally real, we may
choose coordinates so that $H$ vanishes to arbitrarily high order at $x =
0$. These coordinates are not adapted to $M$, so they must be replaced by the
adapted coordinates $(x^1,\dots,x^m,v^1,\dots,v^m,w^{1},\dots,w^{q})$ defined
by
\[
        v^j=  y^j  - H^j(x^1,\dots,x^m)\quad
        w^a = z^{m+a} - H^{m+a}(x^1,\dots,x^m) \,.
\]
Blowup coordinates are then given by the formulas
\[
    v^{\a} =  r p^{\a}, \quad
   v^m =  r,  \quad
   w^{a} = r \zeta = r ( \xi^a + i \eta^a)\,,
\]
where  Greek indices range between $1$ and $m-1$.
Since we only have to  compute $\dc u$ on the set $x^j = 0$, and since $H$ vanishes to high
order, we may assume that $H(x)$ is identically zero in any finite order
computation along $x^j=0$. In particular, up to first order along the set
$x^j=0$, $j=1,\dots,n$, we
have
\[
   z^{\a} = x^{\a} + i r p^{\a}, \quad
   z^m = s + i r,  \quad
   z^{m+a} = r \zeta = r ( \xi^a + i \eta^a) \,.
\]
(To highlight the special role played by the radial parameter $r$, we have written
$z^m = s + ir$.)  A straightforward computation using the chain rule, shows
that
\begin{align*}
\pp^*( \dc u) &= -\left(
               \pd{u}{r} - \frac{p^{\a}}{r} \pd{u}{p^{\a}}
             - \frac{\xi^{a}}{r} \pd{u}{\xi^{a}} - \frac{\eta^{a}}{r} \pd{u}{\eta^{a}}
         \right) d s
 + \left(
         \pd{u}{s} + r p^{\a} \pd{u}{p^{\a}}
          -  \frac{\xi^{a}}{r} \pd{u}{\eta^a} + \frac{\eta^{a}}{r} \pd{u}{\xi^{a}}
    \right) \, d r
\\
& - \frac{1}{r} \,\pd{u}{p^{\a}} \, dx^{\a} + r \pd{u}{x^{\a}}\, dp^{\a}
 - \pd{u}{\eta^{a}}\, d\xi^{a} + \pd{u}{\xi^{a}}\,d\eta^{a} \,.
\end{align*}

We claim that $\tu$ can we written in the form
\begin{equation}
\label{eq:smoothness-1}
          \tu =    r^{2/k} \, U(x',s,p,\xi,\eta,r^{2/k} ) \,,
\end{equation}
where $k>0$ is an integer,  $x'=(x^1,\dots,x^{m-1})$, and $U(x',s, p,\xi,\eta, t)$ is a
differentiable function of $t$ such that
\[
 U(x',s, p,\xi,\eta,0) >0 \,.
\]
Assume this claim for the moment. Then, substituting
\eqref{eq:smoothness-1} into the formula for $\tth = \pp^*\dc u$ and
simplifying gives
\begin{align*}
\tth &=
-r^{ \frac{2}{k}-1}  \left(
                  \frac{2}{k} \, U + r \pd{U}{r} -  p^{\a} \pd{U}{p^{\a}}
                -  \xi^{a} \pd{U}{\xi^{a}} - \eta^{a} \pd{U}{\eta^{a}}
          \right)\, ds   \\
&
+
r^{\frac{2}{k}-1} \left(
                r \pd{U}{s} + r^2  p^{\a} \pd{U}{x^{\a}}
                   +  \eta^{a} \pd{U}{\xi^{a}} - \xi^{a} \pd{U}{\eta^{a}}
           \right)\, dr
 \\
&-r^{ \frac{2}{k}-1} \left(
                \pd{U}{p^{\a}} \, dx^{\a}
          \right)
+ r^{\frac{2}{k}+1 }\left(
                   \pd{U}{x^{\a}} \, dp^{\a}
            \right)
-  r^{\frac{2}{k}}  \left(
             \pd{U}{\eta^{a}}\, d\xi^{a} - \pd{U}{\xi^{a}} \,  d\eta^{a}
          \right)
\end{align*}
But we have already proved that $\tth$ extends smoothly to the set $r=0$
and is nowhere-vanishing. Inspection of the above formula for $\dc u$ shows
that this implies that $k=2$.  Thus, $\tu = r U(x^{\a},p^{\a},s,\xi,\eta,r)$, which is
smooth on all of $\tMC$. The formula $d\tu = U \,dr$ for
$r=0$ shows that $u$ is regular. That that $M$ is totally real and has
dimension $n$ follows from Theorem~\ref{thm:singular}.

It remains only to prove that $\tu$ is of the form \eqref{eq:smoothness-1}.
Because $f$ is real analytic, $\tau$ has a series expansion of the form
\[
           \tau =   f(u) =   a u^k \, (1 + g(u) )
\]
where $a>0$ and $g(u)$ is a smooth function such that $g(0)=0$. Therefore, the equation
\[
          \mysqrt{k}{\tau} = \mysqrt{k}{a}\, u \,\mysqrt{k}{ \left( 1 + g(u)\right)}
\]
can be inverted to show that $u$ is of the form
\begin{equation}
\label{eq:smoothness-2}
             u =    \mysqrt{k}{\tau} \, G( \mysqrt{k}{\tau})
\end{equation}
for $G(t)$ a smooth (in fact, analytic) function satisfying the condition $G(0)>0$.

On the other hand,  $\tau \geq 0$ is smooth and vanishes
precisely on $M$. This, together with the positivity condition
 $d\dc \tau>0$, implies that $\tau$ vanishes precisely to order 2 on
$M$. Therefore, $\tau$ can be expressed in the form
\begin{equation}
\label{eq:smoothness-3}
          \tau =     r^2 \, T(x',s,p,\xi,\eta,r) \,,
\end{equation}
where $T(x',p,\xi,\eta,r)$ is a smooth function and $T(x',p,\xi,\eta,0)>0$.
Combining \eqref{eq:smoothness-2} and \eqref{eq:smoothness-3}, and setting
$U = \sqrt{T}$ results in
the expression~\eqref{eq:smoothness-1}.

Setting $k=2$ in the above formula for $\tth$ and simplifying
yields the identity
\begin{equation}
\label{eq:smoothness-4}
\thS =
  -\left(  U -  p^{\a} \pd{U}{p^{\a}}
                -  \xi^{a} \pd{U}{\xi^{a}} - \eta^{a} \pd{U}{\eta^{a}}
          \right)\, ds
        -         \pd{U}{p^{\a}} \, dx^{\a}  \,.
\end{equation}
At this point, we invoke Theorem~\ref{thm:smoothness} to conclude that
$m=n$ and $q=0$, and that $M$ is totally real.
\end{proof}

\begin{remark}
\label{rem:smoothness-5}
For later reference, we note that because $q=0$,
Formula~\ref{eq:smoothness-4} reduces to the identity
\[
\thS = -\left( U - p^{\a} \pd{U}{p^{\a}} \right)\, ds
       - \pd{U}{p^{\a}} \, dx^{\a}\,.
\]
\end{remark}

\medskip


\section{The Finsler metric on $M$}
\label{sec:finsler}

When $u$ is regular on the blowup, the restriction of
$\tth$ to $\SM$ is a contact form. We now show that that where $u$ is
regular on the blowup, it induces a \emph{Finsler
metric} on $M$ and that the leaves of the \MA\ foliation  intersect $M$
along  geodesics.

\medskip

\subsection{Review of Finsler geometry}
We begin with a quick review of Finsler geometry from the perspective of
contact geometry. For a more complete and more general, exposition of these
ideas the reader should consult the paper of Pang \cite{pang:1990}. Let
$\pi:TM\to M$ denote the tangent space of $M$ and let $T_{0}M \subset TM$
denote the set of non-zero tangent vectors.  A \emph{Finsler metric} on $M$
is a smooth, positive function $\F:T_0 M \to R$ that satisfies the following
two conditions:

\begin{enumerate}
\item[(i)] For all $X \in T_{0}M$ and all $t>0$, $\F(t X) = t \F(X)$.
\medskip
\item[(ii)] The set $S_{p} = \{X \st \F(X)=1 \}$ is strongly convex and
diffeomorphic to a sphere.
\end{enumerate}

Let $x^j$, $j=1,\dots,n$ be local coordinates on $M$ and let
$(x^j,\dot{x}^j)$ be the induced coordinates on the tangent bundle.  The
\emph{Hilbert form} $\thF$ on $T_0M$ is the 1-form defined by the local
formula
\[
       \thF = \pd{F(x,\dot{x})}{\dot{x}^{j}} \, d x^{j} \,,
\]
where the summation conventions are in force.
It is not difficult to show that the convexity condition (ii) is equivalent
to the condition that
\[
      \thF \wedge (d\thF)^{n-1}
\]
be non-vanishing.

The homogeneity of $\F$ implies that $\thF$ is the pullback of a 1-form
on the projective tangent bundle $\SM$, which by abuse of notation we also denote
by $\thF$. To see this, let
\[
       X_R = \dot{x}^{j} \pd{}{\dot{x}^{j}}
\]
denote the radial vector field. We must only show that
\[
X_{R} \i \thF = 0 \text{ and } \Lie{X_{R}}\thF = 0 \,.
\]
The first identity is obvious.  To prove the second, compute as
follows:
\[
\Lie{X_{R}}\thF = X_{R}\i d\thF + d ( X_{R}\i\thF)
   = X_{R}\i d\thF = \dot{x}^{j}  \pdd{F}{\dot{x}^j}{\dot{x}^k} \,dx^{k} =0\,,
\]
where the last equality on the right follows by differentiating
Euler's identity, $\displaystyle \dot{x}^{j} \pd{F}{\dot{x}^{j}} = F$,
with respect to $\dot{x}^k$. We have proved the following lemma:

\begin{lemma}
\label{lem:finsler}
The function $\F$ is a Finsler metric
if and only if the form $\thF$ on $\SM$ is a \emph{ contact form},
i.e.
\[
          \thF \wedge d\thF^{n-1} \neq 0 \,.
\]
\end{lemma}

\begin{remark}
\label{rem:finsler}
The geodesics of a Finsler manifold have an elegant formulation in terms of
the Hilbert form. The \emph{Reeb vector field} of the contact contact
manifold  $(\SM, \thF)$  is the vector field $\XF$ characterized by the
conditions:
\[
      \thF( \XF) = 1 \text{ and } \XF \i d\thF = 0 \,.
\]
The geodesics of $(M,F)$ are the images under the projection map $\pi:\SM
\to M$ of the integral curves of $\XF$. In fact, if $t \mapsto \nu_t(\tp)$
is the integral curve of $\XF$ starting at $\tp \in \SM$, then
$\c: t \mapsto \pi \circ \nu_t(\tp)$ is the unit speed geodesic with
$[\c'(0)] = \tp$.
\end{remark}

\subsection{Construction of the metric}

Let $u$ be a solution of the \MA\ equation and assume that $u$ is regular
on  the blowup of $M$.
Recall that this implies that $M$ is a maximal,
totally real submanifold of $\MC$. Thus the composition of the maps in
\eqref{eq:singular-4} is an isomorphism of vector bundles.

Define $\F:TM_{0} \to \R$ as follows. Let $X$ be a non-zero tangent vector
based at a point $p\in M$. Let $\c(t)$ be curve such that $\c(0)=p$
and $\c'(0) = J X$. Then we set
\begin{equation}
\label{eq:finsler-1}
     \F(X) =   \lim_{t \to 0^{+}}
        \frac{u \circ \c(t)}{t} \,.
\end{equation}
The next proposition shows that $\F$ is a Finsler metric on $M$.

\begin{proposition}
\label{prop:finsler-1}
Suppose that $u$ is regular on the normal blowup of $M$ and let $\thS$
denote contact form on $\SM$ obtained by pulling-back the form $\tth$ to
$\SM$.  Then $\F$ is a Finsler metric and its Hilbert form $\thF$ coincides
with $-\thS$.
\end{proposition}

\begin{proof}
Let $\widetilde{\c}(t)$ be the lift of $\c(t)$ to $\tMC$ defined in
Lemma~\ref{lem:blowup}(i).
Then
\[
      \F(X) =  d\tu(\widetilde{\c}'(0) ) \,.
\]
By (\ref{lem:blowup}(i)), $\F(X)$ depends only on $X$; thus, $\F$ is well defined.
Homogeneity of $\F$ follows from the definition of $\F$.
To see that $F(X)$ is positive, write $\tu$ in the form
\[
          \tu(x,p,r) = r \,U(x,p,r) \,,
\]
where $(x,p,r)$ are blowup coordinates as in \ref{lem:blowup}. Because $u$
is regular on the blowup, $U(x,p,0)$ is strictly positive. Consequently,
\[
            \F(X) = r'(0)\, U\left(x(0),p(0),0\right)
\]
where $\widetilde{\c}(t) = (x(t),p(t),r(t))$.  Finally, observe that
$r'(0)$ is positive because $JX = \c'(0)$ is transverse to $M$.

By Lemma~\ref{lem:finsler}, to conclude the proof we need only show that
$-\thS$ coincides with the Hilbert form of $\F$.  We prove equality via
explicit formulas for both forms using blowup coordinates centered at an
arbitrary point $p_0 \in M$.
Because $p_0$ is arbitrary, we need only
verify equality on the fiber $\pp^{-1}(p_0)$.

Choose  holomorphic coordinates $z^{j} = x^{j} + i\, y^{j}$,
$j=1,\dots,n$, centered at $p_0$ as in the proof of
Theorem~\ref{thm:smoothness}.
In these coordinates, $u$  assumes the form
\begin{equation}
\label{eq:finsler-2}
       u(x,p,r) = r U(x,p,r)
\end{equation}
and, by Remark~\ref{rem:smoothness-5}, the form $\thS$ assumes the form
\begin{equation}
\label{eq:finsler-3}
  \thS =
   -\left(U  - p^{\a} \pd{U}{p^{\a}} \right)ds -  \pd{U}{p^{\a}} dx^{\a}
  \,.
\end{equation}

We next focus on the computation of $F$ and $\thF$. Let
\[
     X  = \dot{x}^{1}\, \pd{}{x^1} + \dots + \dot{x}^{n-1}\, \pd{}{x^{n-1}}
             + \dot{s} \,\pd{}{s}
\]
denote a tangent vector to $M$ at $p_0$. If  the
ray generated by $JX$ is in the coordinate patch of $\tMC$, then
$\dot{s}>0$.   Then by \eqref{eq:finsler-2}
and Lemma~\ref{lem:blowup},
\[
        F(X) =    \dot{s} \,U\left( x,p,0 \right) \,,
\]
where $p^{\a} = \dot{x}^{\a}/\dot{s}$.  The Hilbert form of $F$ is
therefore given by
\[
   \thF =     \pd{F}{\dot{x}^j} \, d x^j
        = \left(U  - p^{\a} \pd{U}{p^{\a}} \right) ds + \pd{U}{p^{\a}}
        dx^{\a} \,.
\]
Comparing  this formula with \eqref{eq:finsler-3} yields the equality
 $\thF = - \thS$ and concludes the proof of the proposition. \end{proof}

\begin{proof}[Proof of Theorem~\ref{thm:main}(a)]
The proof is a corollary to Proposition~\ref{prop:finsler-1}. Because the forms
$-\thS$ and $\thF$ coincide, the Reeb vector field of $\thF$ coincides with
the restriction to $\SM$ of the vector field $X$ defined by
Equation~\eqref{eq:singular-2}. But Lemma~\ref{lem:singular} shows that the
projection onto $M$ of the integral curves of $X$ are the intersections of
leaves of the \MA\ foliation with $M$.\end{proof}


\section{Construction of regular \MA\ models}
\label{sec:models}

In this section, we prove parts (b) and (c) of Theorem~\ref{thm:main}.
Our proof is a generalization of a construction of Lempert-Sz\"{o}ke
\cite{lempert-szoke:1991a}.

Before beginning the proof, we make a few preliminary observations.  Recall
that, because $M$ is totally real, the complex structure tensor $J$ induces
an analytic isomorphism between the projective tangent bundle of $M$ and
the projective normal bundle $\S Q$, which is, by construction, the
boundary of $\tMC$.  We may, therefore, identify the boundary $\SM$ of
$\tMC$ with the projective tangent bundle of $M$.

Thus far, we have worked in the smooth category; we now introduce the
further assumption that all data are real analytic.  Specifically, let
$(M,\F)$ denote a compact, real analytic manifold with a real analytic
Finsler metric. Then the oriented projective tangent bundle $\SM$ is also
real analytic, as are the Hilbert form $\thF$ and the Reeb vector field
$\XF$. If follows that the flow of $\XF$,
\[
     \nu: \SM \times \R \to \SM \st (p,t) \mapsto \nu_t(p)
\]
defines a real analytic family of diffeomorphisms of $\SM$.

Next let $\MC$ denote the complexification of $M$.  By construction, $M$ is
an analytic, $n$-dimensional, totally real submanifold of its
complexification $\MC$, and  any real analytic atlas for $M$ extends to define a
holomorphic atlas for $\MC$. Using this atlas to define the normal blowup
as in Section~\ref{sec:blowup} immediately shows that $\tMC$ has real
analytic boundary and that the blowdown map
\[
                 \pp:\tMC\to\MC
\]
is real analytic.
With this identification, we have the following diagram of
real analytic maps:
\begin{equation}
\label{eq:models-1}
     \SM \times \R \mapright{\nu} \SM
                   \hookrightarrow \tMC
                   \mapright{\pp} \MC  \,.
\end{equation}
We are now going to  extend this map to the domain $\SM \times \C$ by analytic
continuation and use the extension to define a solution $u$ of the \MA\
equation.
The map \eqref{eq:models-1} gives a real-analytic family of curves,
$\gamma_{\tp}$, $\tp \in \SM$, defined by
\begin{equation}
\label{eq:models-2}
   \gamma_{\tp}: \R \to M \subset \MC \st t \mapsto \pp \left(\nu_{t}(\tp) \right)\,,
\end{equation}
and each curve is both a geodesic in the Finsler manifold $(M,\F)$ and a
real analytic curve in $\MC$.  By virtue of the second property, each of
these curves can be holomorphically extended to a holomorphic curve defined
on a neighborhood of of $\R$ in $\C$. The next lemma shows that the
extension is  uniform over all of $\SM$.

\begin{lemma}
\label{lem:models} There exists a real number $R>0$ and a real analytic
extension
\[
   \nu^{\C} : \SM \times H_{R} \to \tMC
\]
of $\nu$, where $H_{R} = \{ s + i r \st 0 \leq y < R \}$. The map
$\nu^{\C}$ has the following properties:

\begin{enumerate}

\item[(i)] For each $p \in \SM$, the map $z \mapsto \pp\circ \nu_{C}(p,z)$ is a
holomorphic immersion.

\item[(ii)] The map $\mu: \SM \times [0,R) \to \tMC$ defined by the formula
\[
          \mu(p,r)   = \nu^{C}(p, 0 + r i)
\]
is a real analytic diffeomorphism onto its image.
\end{enumerate}
\end{lemma}

\begin{proof}
Choose a point $p \in \SM$. Because $\c_{p}$ is real analytic, for
sufficiently small $\epsilon>0$, it has a holomorphic extension
$\c^{\C}_{p}:V \to \MC$, where $V_{\epsilon} = \{ z = s + i r \st
|s|<\epsilon,\, 0 \leq r<\epsilon$. It is easy to check that
$\c^{\C}_{p}$ lifts to a real analytic map
\[
       \tc^{\C}_{p}:V \to \tMC
\]
which is an extension $\nu$. By analytic dependence of
$\c_{p}$ on $p$ and compactness of $\SM$, there exists a real
number $R>0$ such that $\c_{p}^{\C}$ is defined on $V_{R}$ for all
$p \in \SM$. We now have a real analytic map
\[
      \nu^{\C}\st \SM \times V_{R} \to \tMC \,,
\]
which is holomorphic in the second factor. The one-parameter
identity $\nu_{t+s} =\nu_t \circ \nu_s$ then allows us to extend the map to
all of $\SM \times H_R$ as the composition
\[
     \nu^{\C}_{s+i r}(p) = \nu^{\C}_{ s/k +i r} \circ \dots \circ
                     \nu^{\C}_{ s/k +i r}(p) \,,
\]
where the integer $k$ is chosen so that $|s/k| < R$.

Property~(i) of $\nu$ follows by construction. To prove property (ii),
first observe that $\mu$ is the identity map on $\SM \times \{0\}$. We,
therefore, need only show that the derivative of $\mu$ is injective on all
of $\SM$. It then follows (after shrinking $R$ if necessary) that $\mu$ is
a diffeomorphism, as claimed.  But because $\mu$ is the identity on $\SM$,
it follows that $\mu_{*}$ is injective if and only if the vector field
$\mu_{*}(\d/\d r)$ is transverse to $\SM$. It suffices to show that the
projection $\pp{}_{*}\mu_{*}(\d/\d r)$ is transverse to $M$. But this is
clear, for by construction
\[
    (\pp{} \circ\mu)_{*}(\d/\d r) = \frac{d}{dr} \gamma^{\C}_{p}(i r) =
   J \nu'_{t}(p) \,.
\]
This completes the proof of the lemma.\end{proof}

Replace $\MC$ by the image of $\mu$, and let
$\tu:\tMC\to\R$ be the smooth function defined by the formula
\begin{equation}
\label{eq:models-3}
     \tu:  \tMC \mapright{\,\mu^{-1}} \SM \times [0,R) \mapright{\pi_2}
     \R,
\end{equation}
where $\pi_2$ is projection onto the second factor. Because $\tu$ vanishes
on $\SM$, it descends to a continuous function $u$ on $\MC$.
To complete the proof of  Theorem~\ref{thm:main}(b), we need only show
that  $(\MC,M,u)$ is a real analytic \MA\ model for $(M, \F)$.
We need only  check that the following  conditions are satisfied (after
possibly further shrinking  $R$) :

\begin{enumerate}

\item[(i)] $\tu$ is smooth on all of $\tMC$;

\item[(ii)] $u$ induces the Finsler metric $\F$. Specifically, choose a
tangent vector $X$ and let $\tc(t)$ be the lift to $\tMC$ of a smooth curve
$\c(t)$ with $\c'(0) = J X$, then $\F(X) =   d\tu(\tc'(0))$;

\item[(iii)] $d\dc u^2>0$ on $\MC \setminus M$;

\item[(iv)] $(d\dc u)^n = 0$ on $\MC \setminus M$;

\item[(v)] $\th = \dc u$ lifts  to a smooth form $\tth$ which extends
smoothly to  all of $\tMC$
and which satisfies the inequality $d\tu \wedge \tth \wedge(d\tth )^{n-1} \neq 0$.

\end{enumerate}

Properties (i) and (ii) follow immediately from the constructions above.

\smallskip

To verify condition (iii), choose an arbitrary point  $p\in\SM$ and choose
blowup  coordinates $(x^{\a},s,p^{\a},r)$ centered at $p$
with $z^{\a} = x^{\a}+i r p^{\a}$, $z^{n} = s + i r$ holomorphic
coordinates on $\MC$. We claim that
complex Hessian of
\[
  H_{\C}(u^2) =
\begin{pmatrix}
  \pdd{u^2}{z^{\a}}{z^{\bb}} & \pdd{u^2}{z^{\a}}{z^{\nb}} \\
  \pdd{u^2}{z^n}{z^{\bb}} & \pdd{u^2}{z^n}{z^{\nb}}
   \end{pmatrix}
\]
extends continuously to a positive definite matrix on a neighborhood of
$p$. Assume the claim for the moment. By compactness of $\SM$, there
is an open neighborhood $U\subset \MC$ of $M$ such that
$d\dc u^2>0$ on $U \setminus M$. By  shrinking $R$ if necessary, we may
assume that $U = \MC$.

To prove the claim, observe that by construction
 $u= r U(s,r,x^{\a},p^{\a})$, where $U$ is smooth and
 $U(s,0,x^{\a},p^{\a})>0$.  Noting that
\[
x^{\a} = \frac{1}{2}( z^{\a} + z^{\ab}) \,,\quad
     p^{\a} = \frac{ z^{\a} - z^{\ab}}{z^{n} - z^{\nb}} \,,\quad
     s = \frac{1}{2} (z^n + z^{\nb})\,,\quad
     r = \frac{1}{2i}(z^n - z^{\nb})
\]
and applying the chain rule to $u^2$ yields the formulae
\begin{equation}
\label{eq:HC}
\begin{cases}
\displaystyle
\pdd{u^2}{z^{\a}}{z^{\bb}} &=
\displaystyle
 \frac{1}{4}\pdd{(U^2)}{p^{\a}}{p^{\b}}
\\
\displaystyle
\pdd{u^2}{z^{\a}}{z^{\nb}} &=
\displaystyle
\frac{1}{4}\pd{}{p^{\a}}
        \left( 2 U^2 -  p^{\b} \pd{(U^2)}{p^{\b}}\right)
\\
\displaystyle
\pdd{u^2}{z^{n}}{z^{\nb}} &=
\displaystyle
\frac{1}{4} \left(
 2 U^2 - 3 p^{\a} \pd{(U^2)}{p^{\a}} +  p^{\a} p^{\b}
 \pdd{(U^2)}{p^{\a}}{p^{\b}}
\right)
\end{cases}
\end{equation}
for $r=0$.

\medskip

On the other hand, (ii) implies that the Finsler metric $F$ has the form
\[
           F(x,\dot{x}) = \dot{x}^n U(s,0, x^{\a}, \dot{x}^{\a}/\dot{x}^n) \,.
\]
The convexity for $F$ implies that the real Hessian
\[
         H_{\R} (F^2) =
 \begin{pmatrix}
\pdd{F}{\dot{x}^{\a}}{\dot{x}^{\b}} & \pdd{F}{\dot{x}^{\a}}{\dot{x}^{n}}\\
\pdd{F}{\dot{x}^n}{\dot{x}^{\a}}    & \pdd{F}{\dot{x}^n}{\dot{x}^{n}}
 \end{pmatrix}
\]
is positive definite for all $\dot{x} \neq 0$.  A straightforward
computation shows that
\begin{equation}
\label{eq:HR}
\begin{cases}
\displaystyle
\pdd{(F^2)}{\dot{x}^{\a}}{\dot{x}^{\b}} &=
\displaystyle
\pdd{(U^2)}{p^{\a}}{p^{\b}}
\\
\displaystyle
\pdd{(F^2)}{\dot{x}^{\a}}{\dot{x}^{n}} &=
\displaystyle
\pd{}{p^{\a}}\left( 2 U^2 - p^{\b} \pd{(U^2)}{p^{\b}} \right)
\\
\displaystyle
\pdd{(F^2)}{\dot{x}^{n}}{\dot{x}^{n}} &=
\displaystyle
 2 U^2  - 3 p^{\a} \pd{(U^2)}{p^{\a}}
                    + p^{\a} p^{\b} \pdd{(U^2)}{p^{\a}}{p^{\b}}
\end{cases}
\end{equation}
Comparison of \eqref{eq:HC} and \eqref{eq:HR} shows that
$H_{\C}(u^2) = \frac{1}{4} H_{\R}(F^2)$ at $r=0$. Consequently, the complex Hessian of $u^2$ is
positive definite in a neighborhood of $p \in \SM$.

\smallskip

To verify condition (iv), recall that by Lemma~\ref{lem:models} every
point of $\MC \setminus M$, is contained in the image of a holomorphic
curve of the form $z \mapsto \nu^{\C}(p,z)$.  By definition, $u \circ
\nu^{\C}(p,z) = \Im(z)$, showing that the pull-back of $d\dc u$ to the
curve vanishes. Together with (iii), this shows that
 $d\dc u$ has rank strictly less than $n$.

\smallskip

To verify condition (v), we  first show that the form $\tth =
\pp^{*}\dc u$ extends to all of $\tMC$. To see this note that by
construction, the vector field on $\tMC$
$
           Y  =   \mu_{*} \pd{}{r}
$
satisfies the identity
$
          Y \i \tth =   0
$
on the set $\tMC \setminus \SM$, and the computation
\[
        \Lie{Y}\tth = Y \i d\tth = 0
\]
shows that the Lie derivative vanishes.  Consequently,
we need only verify (iv) on $\SM$.  But since $Y(\tu)=1$, we need only
show that $\thS$, the pull back of $\tth$ to $\SM$,  is a contact form.
But the computations leading to the formula~\eqref{eq:finsler-3} all apply here,
showing that $\thF = -\thS$. Non-degeneracy follows, concluding the proof of
part (b) of Theorem~\ref{thm:main}.

\medskip

To prove part (c) of Theorem~\ref{thm:main}, first suppose that $\Phi:\MC \to
\MC$ is a biholomorphism between two \MA\ models $(\MC,M,u)$ and $(\MC'
,M',u')$ such that $u = u' \circ \Phi$. Equation~\eqref{eq:finsler-1} shows
the $\Phi$ restricts to an isometry between $(M,\F)$ and $(M',\F')$.
Conversely, any analytic isometry between real analytic Finsler manifolds
 $(M,\F)$ and $(M',\F')$ extends uniquely to a biholomorphism $\Phi:\MC \to
 \MC'$ between their complexifications.

Therefore, we need only show that  $u = u' \circ \Phi$, which we can do by
proving equality  on each leaf of $\Fol$. To this end,
let $\c^{\C}: H_R \to \MC$ be the holomorphic parameterization of a leaf
given above. Then $U(z) = u \circ \c^{\C}(z)$
and $U'(z) = u'\circ\Phi\circ \c^{\C}(z)$ are both real
analytic solutions of the  initial
value problem
\[
    \pd{^2U}{r^2} = - \pd{^2U}{s^2} \,,\quad
      U(s,0) = 0\,,\quad
      \pd{U(s,0)}{r} = 1 \,.
\]
By the Cauchy-Kovaleskaya Theorem, it follows that $U(s,r)=r$. Uniqueness
follows,  completing the proof of Theorem~\ref{thm:main}.




\end{document}